\documentstyle {article}

\newcommand{\be}{\begin{equation}}
\newcommand{\ee}{\end{equation}}
\newcommand{\zz} {\overline{z}}
\newcommand{\nn} {\noindent}
\begin{document}
\title{  Deformed harmonic oscillators : coherent states and Bargmann representations}
\author{Mich\`{e}le IRAC-ASTAUD and Guy RIDEAU \thanks{Presented at the 6th Colloquium "Quantum Groups and Integrables Systems", Prague, 19-21 June 1997}\\ 
Laboratoire de Physique Th\'{e}orique de la mati\`{e}re condens\'ee\\
Universit\'{e} Paris VII\\2 place Jussieu F-75251 Paris Cedex 05, FRANCE}
\date{}
\maketitle

\begin{abstract}
Generalizing the case of the usual harmonic oscillator, we look for Bargmann representations corresponding to deformed harmonic oscillators. Deformed harmonic oscillator algebras are generated by  four operators $a, a^\dagger, N$ and the unity $1$ such as  $[a,N] = a ,\quad [a^\dagger ,N] = -a^\dagger$, $a^\dagger a = \psi(N)$ and $\quad aa^\dagger =\psi(N+1)$. We discuss the conditions of existence of a scalar product expressed with a true integral on the  space spanned by the eigenstates of $a$ (or $a^\dagger$). We give various examples, in particular we consider functions $\psi$ that are linear combinations of $q^N$, $q^{-N}$ and unity and that correspond to q-oscillators with Fock-representations or with non-Fock-representations. 
\end{abstract}

\section{Introduction}

\nn The deformed harmonic oscillator algebra is defined by  four operators : the annihilation operator $a$, the creation operator $a^\dagger$ , the energy operator $N$ and the unity $1$  satisfying the following commutation relations :

\begin{equation}
[a,N] = a ,\quad [a^\dagger ,N] =
-a^\dagger ,\quad
a^\dagger a = \psi(N) , \quad aa^\dagger =\psi(N+1)
\label{a3}
\ee

\nn where $\psi$ is a real analytical function.

\nn Let us give some examples :

\nn - the usual harmonic oscillator corresponds to $
\psi_{usual}(N) = N +\sigma$.

\nn - the q-oscillator  defined by $aa^\dagger - q a^\dagger a = q^{-N}$ corresponds to 

\be
 \psi_{qosc} (N) = - q^{-N}/(q-q^{-1}) + \sigma q^N /(q-q^{-1}), \forall \sigma. 
\label{qosc}
\ee

\nn - the q-oscillator  defined by $aa^\dagger - q a^\dagger a = 1$ corresponds to 

\be
 \psi ^\prime_{qosc} (N) = (1-q)^{-1} + \sigma q^N , \forall \sigma.
\label{qosc2}
\ee

\nn  Generalizing the pioneer work of Bargmann \cite{Bargmann} for the usual harmonic oscillator, the purpose of this paper is to study if  the deformed harmonic oscillator defined by (\ref{a3}) admits representations on one space of complex variable functions.  
The scalar product of the representations we are looking for, is written with a true integral  contrarily to many studies where a q-integration occurs. 

\nn In section 2, we describe  the irreducible representations on the basis of the eigenvectors of $N$, they depend on the zeros of $\psi$. We then discuss the existence of the coherent states that are defined as the eigenstates of the operators $a$ (or $a^\dagger$). 
We study the possibility of  Bargmann representations, in section 3. We show on various examples how the construction works  :  section 4 is devoted  to strictly positive function $\psi$ , other cases are considered in section 5.

\section{Coherent States}

\nn Let us recall that the representations on a basis of eigenvectors of $N$ are characterized by the zeros of the function $\psi$
\cite {cohbarg} :

\nn 1){\bf $\psi$ has no zero.}

\nn The inequivalent representations are labelled by the decimal part of $\mu$ and are defined by :

\be
\left\{
\begin{array} {ll}
a^\dagger \mid n >  =& (\psi (\mu +n+1))^{1/2} \mid n+1 > \\ 
a\mid n >  =& (\psi (\mu + n))^{1/2} \mid n-1 > ,\quad n \in Z \\ 
N \mid n >  =&  (\mu +n )\quad \mid n > 
 \end{array}
\right.
\ee

\nn The spectrum of $N$, $Sp N$, is $\mu + Z$. The operator $N$ has
 no lowest and no highest eigenstates. These representations, thus, are non equivalent to  Fock-representations. 

\nn {\bf 2) $\psi$ has zeros.}

\nn We are interested in the intervals where $\psi$ is positive :

\nn {\bf a) finite intervals }

\nn We can associate one representation to the intervals that have a length equal to an integer. Let us denote $\nu_-$ and $\nu_+$ two integers such as $\psi (\mu +\nu_-) =0$ and $\psi (\mu +\nu_+ +1) =0$. If $\psi (x) >0$ when $x \in ]\mu +\nu_-,\mu +\nu_+ ]$,
the spectrum of $N$ is $[ \mu + \nu_-,\mu+\nu_+] \bigcap Z +\mu$.

\nn {\bf b) infinite intervals}

\nn The representations are similar to the  Fock-representation of the usual harmonic oscillator :

\nn - $Sp N = \mu + \nu_- + N^+$ if $\psi(\mu +\nu_-)=0$ and $\psi(x) > 0 $ when $x > \mu +\nu_-$  

\nn - $Sp N =\mu + \nu_+ + N^-$ if $\psi(\mu +\nu_+ +1)=0$ and $\psi(x) > 0 $ when $x \leq \mu +\nu_+ $.

\nn The first step to build a Bargmann representation requires to study the coherent vectors.
We call {\bf coherent states} $\mid z>$ , the eigenvectors of the operator $a$ (or $a^\dagger$) with eigenvalues $z$. 

\nn - When the spectrum of $N$ is upper bounded,  $a$ has no eigenvectors and when the spectrum of $N$ is lower bounded,  $a^\dagger$ has no eigenvectors.
Therefore, in the case (2.a)  as the spectrum of $N$ is finite, $a$ and $a^\dagger$ have no eigenvectors , hence  no Bargmann representation exists.

\nn - When the spectrum of $N$ is no upper bounded,
the  eigenvectors $\mid z>$ of $a$ take the form~:

\be
\left\{
\begin{array}{lll}
\mid z>=&\sum_{n=-1}^{\nu_-}z^{n}(\psi(\mu +n)!)^{1/2}\mid n>+\sum_{n=0}^{\infty}\frac{z^n}{ \psi(\mu +n)!^{1/2}}\mid n>,& \nu_- < 0\\
&&\\
\mid z>=&\sum_{n=\nu_-}^{\infty}z^n (\psi(\mu +n)!)^{-1/2}\mid n>,&\nu_- \geq 0\\
\end{array}
\right.
\label{z}
\ee

\nn with the convention $\psi(\mu )! =1$ and $\nu_- = -\infty$ if $Sp N = \mu +Z$.

\nn The domain $D$ of existence of the coherent states depends on the function $\psi$. Indeed, $\mid z>$ belongs to the Hilbert space $\cal {H}$ spanned by the basis $\mid n>$ only if the series in the right hand side of (\ref{z}) are convergent in norm. 

\nn The eigenvectors of $a$ exist if :

\nn - $Sp N = \mu+Z$ and $ \psi(+\infty)>  \psi(-\infty)$, the domain  $D$ of existence of the coherent states is $D = \{ z; \psi(-\infty) < \mid z \mid ^2< \psi(+\infty)\}$ or 

\nn - $Sp N = \mu +\nu_- + N^+$, then $D = \{ z;\mid z \mid ^2 < \psi(+\infty)\}$ . 

\nn The eigenvectors of $a^\dagger$ exist if :

\nn - $Sp N = \mu +\nu_+ + N^-$, then $D = \{ z;\mid z \mid ^2 < \psi(-\infty) \}$ or

\nn - $Sp N = \mu+Z$ and $\psi(+\infty) < \psi(-\infty)$, then  $D = \{ z; \psi(+\infty)< \mid z \mid ^2 < \psi(-\infty) \}$.

\nn The part plays by $a$ and $a^\dagger$ being analogous, in the following we restrict to the case where the eigenstates of $a$ exist.

\nn Although $\mu$ is a significant quantity as labelling inequivalent representations, it does not play a  part in the present problem. 
So we simplify the notation in assuming $\mu =0$ from now on. 

\section{ Bargmann representation}

\nn Following the construction \cite{Bargmann}, in the Bargmann representation any state $\mid f>$ of $\cal{H}$
is represented as the function of a complex variable $z$, $f(z) = < \zz \mid f >$, with a Laurent expansion   
 on the domain $D$ of definition of the eigenvectors of $a$. The space of the Bargmann representation $\cal{S}$ is constituted with holomorphic functions in $D$, strongly depending on $\psi$. 

\nn A  Bargmann representation exists if  we can exhibit  a positive real function $F(x)$ such as 

\be
\int F(z \zz ) \mid \zz ><\zz \mid dz d\zz =1
\label{1}
\ee

\nn where the integration is extended to the whole complex plane and where $F(\mid z \mid ^2 )$ contains the characteristic function of the domain $D$ of existence of the  coherent states.

 \nn We easily prove  that, in  this representation, $a^\dagger$  is the multiplication by $z$ and $N$ is the operator $zd/dz$, as in the usual case, while  $a$ is the operator $z^{-1}\psi(zd/dz)$. Let $G(\zeta z)$ be equal to $<\zz \mid \zeta >$, it is the function of  $\cal{S}$ corresponding to the coherent state $\mid \zeta >$. If (\ref{1}) holds, it appears to be a reproducing kernel. Moreover, it fulfills the following equation
\be
x G(x) = \psi(x \frac{d}{dx}) G(x)
\label{GG}
\ee

\nn Let us introduce the Mellin transform $\hat{F}(\rho)$ of the weight function $F(x)$ :

\be
\hat {F}(\rho)=\int_0^{\infty} F(x) x^{\rho -1} dx = \int_{D } F(x) x^{\rho -1} dx
\label{mel}
\ee

\nn We easily deduce that $F(x)$ must verify the following condition~:

\be
\hat{F}(n+1)= \left\{ \begin{array}{lll}
 \psi(n)!, & n\geq 0&\\
 (\psi(n)!)^{-1}, & n <0& , n \in Sp N
\end{array}
\right.
\label{moments}
\ee

\nn Let us remark that 

\be
\hat {F}(\rho) \leq \hat {F}(n)+\hat {F}(n+1), \quad n \leq Re \rho <n+1
\label{nn}
\ee

\nn as $F(x)$ is a positive function.
We see that, the Mellin transform of $F$ exists for the integers belonging to the spectrum of $N$ and therefore, due to (\ref{nn}), for any $\rho$  such as $Re \rho$ is greater than the lowest bound $\nu_-$ of $Sp N$. 

\nn Formula (\ref{moments}) is equivalent to

\be
\hat{F}(n+1)=\psi (n) \hat{F}(n),\mbox{with} \quad \hat{F}(1)=1
\label{M}
\ee

\nn which ensures that the operators $a^\dagger = z$,  $a = z^{-1}\psi(zd/dz)$ be adjoint on the basis $\mid n>$ of ${\cal H}$.
 The function $\psi$ being given,  $\hat{F}$ verifying (\ref{M}) must be interpolated in order to get $F$ as the Mellin inverse of (\ref{mel}). 

\nn It can be proved \cite{cohbarg} that the only interpolation of (\ref{M}) to be considered  is  the simplest one :

\be
\hat {F}(\rho+1)= \psi(\rho)\hat {F}(\rho)
\label{mel1}
\ee

\nn We assert that a Bargmann representation can be defined on a deformed harmonic oscillator algebra  if the coherent states exist and if  at least one solution of (\ref{mel1}) has  a positive Mellin inverse on D.

\nn Furthermore, let us remark that (\ref{mel1}) can be written :

\be
\int F(x) x^{\rho} dx = \int F(x) \psi(x\partial _x +1) x^{\rho -1} dx
\ee

\nn The right member is equal to $\int \left( \psi(-x\partial_x) F(x)\right) x^{\rho -1} dx$ if the integration by parts (or the change of variables) can be done without extra terms and this gives :

\be
x F(x) = \psi(-x \frac{d}{dx}) F(x)
\label{FF}
\ee

\nn This equation, when it holds, can be used to study the positivity of $F(x)$, when we cannot obtain an explicit expression for this function. 

\nn The general procedure developed in this section is illustrated by some examples in the next sections.

\section {$\psi$ strictly positive }

\subsection{$\psi(N) = \lambda q^{-N}, q \leq 1 , \lambda >0$.}

\nn  When $\lambda = (q^{-1}-q)$, $\psi(N) = \lambda q^{-N}$ corresponds to the only q-oscillator algebra (\ref{qosc}) with non-Fock representations and with coherent states. We easily verify that $Sp N = Z$ and $D= C -\{0\}$. 

\nn In this case, the functional equations (\ref{mel1}) and (\ref{M}) are equivalent.

\nn We obtain the particular solution of (\ref{mel1}):

\be
\hat{F}(\rho)= \exp\left( \rho \ln\lambda -\frac{1}{2}(\rho^2-\rho)\ln q \right)
\ee

\nn  which has the Mellin inverse  \cite{kowalski} \cite{canada}:

\be
F(x) = \exp\left(\frac{\ln^2(x\lambda ^{-1})}{2\ln (q)}-\frac{\ln(x\lambda ^{-1})}{2}\right)
\label{F0}
\ee

\nn Following  Bargmann's procedure  \cite{Bargmann}, we can prove the closedness of the operators $z$ and $z^{-1}q^z\frac{d}{dz}$.  Moreover, we establish that these operators have the same domain of definition and are  mutually adjoint.

\subsection{ $\psi (x) = \exp \left( \sum_0 ^{2p+1} a_n x^n \right), \quad a_{2p+1}>0$.}

\nn A solution of Equation (\ref{mel1}) in this case , is given by :

\be
\hat{F}(\rho ) = \exp \left( \sum_0 ^{2p+1} \frac{a_n}{n+1}B_{n+1}(\rho)
  \right) 
\label{B}
\ee

\nn where $B_{n+1}(\rho)$ are the Bernouilli polynomials. 

\nn We can prove that the Mellin inverse Mellin  of $\hat{F}(\rho)$ exists only if $p$ is an even number. The function $F(x)$ thus obtained is always real, but not necessarily positive. Nevertheless, in specific cases, for example when the exponent in (\ref{B}) contains only the term of highest degree, $F(x)$ is actually strictly positive.
Since $ \psi(n+1)/\psi(n)$ grows indefinitely as $n \rightarrow \pm \infty$,  the domain of $z^{-1}\psi ( z \frac {d}{dz})$ is included in the domain of $z$ but cannot be identical. 

\nn It is worthwhile to underline that, for the functions $\psi$ considered in this subsection and corresponding to odd $p$, we have obtained  examples of deformed oscillator algebras which admit coherent states but for which  Bargmann representations do not exist.

\subsection { $\psi (x) = a + q^x , \quad q>1,\quad a> 0$.}

\nn This function $\psi$ is mainly involved in the study of the q-oscillator (\ref{qosc2}) with non-Fock representations. The domain of existence of the coherent states is a ring on the complex plane $a \leq \mid z \mid^2$.
  We have as a convenient particular solution of (\ref{mel1}),  the following entire function :

\be
\hat{F}(\rho)=  a^{\rho} \prod_{p \geq 0}(1+ a^{-1} q^{\rho-p-1})
\label{aa}
\ee

\nn but its  Mellin inverse is not  a true function. Writing (\ref{aa}) on the form :

\be
\hat{F}(\rho) = a^{\rho}\left(1+\sum_{n\geq1} \frac{a^{-n}q^{n\rho}}{(q-1) \cdots (q^n-1)}\right)
\ee

\nn we verify that it is the Mellin transform of the measure :

\be
F(x)= \sum_{n\geq 0} \frac{a^{-n}}{(q-1) \cdots (q^n-1)} \delta (\ln a+\ln q -\ln x)
\ee

\nn Therefore, in this case we obtain a Bargmann representation if we accept the weight function to be a true measure.

\section {$\psi$ vanishes } 

\subsection{$\psi (x) = [x] \equiv (q^x -q^{-x})/(q-q^{-1}), \quad q\in R$.}

\nn This case corresponds to the q-oscillator defined by (\ref{qosc}) with $\sigma =1$. We verify that $Sp N = N^+$ and $D$ is the whole complex plane.
 The resolution of the identity was obtained with the q-integration \cite{gray}, the weight function being $Exp_q(-x))=\sum_{n\geq 0}(-x)^n/[n]!$.

\nn We here look for a Bargmann representation where the scalar product involves a true integral. Equation (\ref{FF}) and (\ref{mel1}) are equivalent in this case and $Exp_q(-x)$ is a particular solution of (\ref{FF}). It is not a suitable weight function because it vanishes when $x >0$ and we have to look for another solution. As the problem is symmetric under the change $q$ into $q^{-1}$, we now assume without loss of generality that $q >1$.

\nn Let us write a solution of (\ref{mel1}), $\hat{F}$ on the form :

\be
\hat{F}(\rho)= \phi q^{\frac{\rho}{2}(\rho -1)}(q-q^{-1})^{-\rho} \hat{f}(\rho)
\label{fr}
\ee

\nn The function $\hat{f}(\rho)$ must verify

\be
\hat{f}(\rho+1)=(1-q^{-2\rho})\hat{f}(\rho)
\ee

\nn and is given by :

\be
\hat{f}(\rho)= \sum_{n\geq 0}\frac{q^{-2n\rho}}{(1-q^{-2}) \cdots (1-q^{-2n})}
\label{fr1}
\ee

\nn The condition $\hat{F}(1)=1$, furnishes the normalization factor :

\be
\phi = (q-q^{-1})\left( 1+\sum_{n>0}\frac{q^{-2n}}{(1-q^{-2}) \cdots (1-q^{-2n})}\right)^{-1}
\label{fr2}
\ee

\nn Putting (\ref{fr1}) and (\ref{fr2}) in (\ref{fr}), we obtain $\hat{F}(\rho)$ 
, and then we can calculate its Mellin inverse :

\be
F(x)= \frac{\exp(-\frac{1}{2\ln q}(\ln x +\ln(q-q^{-1})+\frac{1}{2}\ln q)^2)}{ 1+\sum_{n>0}\frac{q^{-2n}}{(1-q^{-2}) \cdots (1-q^{-2n})}}\sum_{n\geq 0}\frac{q^{-n(2n+1)}((q-q^{-1})x)^{-2n}}{(1-q^{-2}) \cdots (1-q^{-2n})}
\label{Fx}
\ee

\nn This function is positive and thus can be adopted as  weight function of a Bargmann representation. Therefore, in this case two resolutions of the unity coexist one involving a q-integration and the weight function $Exp_q(-x)$ and another written with a true integral and the weight function $F(x)$ given in (\ref{Fx}). Let us remark that $F(x)$  (\ref{Fx}) verifies the same functional equation (\ref{FF}) as  $Exp_q(-x)$ and thus this function can be the weight function of  one resolution of the identity involving the q-integration. Due to the positivity of $F(x)$, the q-integration is  performed on the whole positive axis, contrarily to  the case where the weight function is $Exp_q(-x)$ \cite{gray}. 

\nn The same is true for the q-oscillator (\ref{qosc2}) with $\sigma = (q-1)^{-1}$.

\subsection {$\psi (x) =(x) \equiv (q^x-1)/(q-1)$, with $ q > 1$.}

\nn In this case, we can prove that the new exponential ${\cal E}xp_q(-x))=\sum_{n\geq 0}(-x)^n /(n)! $ is always positive. Therefore, coexist two resolutions of the identity ,one involving a true integral and the weight function $F(x) = ({\cal E}xp_q(qx))^{-1}$ and one with the new q-integration, the weight function being ${\cal E}xp_q(-x)$. 

\section{Conclusion}

\nn We have studied the possibility of Bargmann representations for any deformed oscillator algebra characterized by a function $\psi$. We gave the conditions to be verified by this function for admitting representations with coherent states. We get the unique functional equation to be satisfied by the Mellin transform of the weight function defining the scalar product. We were able to get definite and positive answer in many cases including in particular some types of q-oscillators. Although we did not succeed in obtaining a general characterization of the function $\psi$ leading to Bargmann representations, we underline two points~: 

\nn - We exhibited cases where the Bargmann representations do not exist even when coherent states do ( subsection (4.2));

\nn - The analysis of subsection (4.3) showed that the scope of our study have to be extended up to include true measures for writing the scalar product.

\nn Finally let us stress that we have obtained  scalar products for the Bargmann representations of the usual q-oscillators, involving  true integrals instead of  q-integrations as previously proposed in literature.

\end{document}